\documentstyle[12pt]{article}
\begin{document}
%
\newtheorem{Satz}{Theorem}
\newtheorem{Theorem}[Satz]{Theorem}
\newtheorem{Lemma}[Satz]{Lemma}
\newtheorem{Kor}[Satz]{Corollary}
\newcommand{\nc}{\newcommand}
\nc{\abs}{\vskip 8pt}
\nc{\bbs}{\vskip 6pt}
\nc{\Seite}{\newpage} 
\nc{\nl}{\newline}
\nc{\noi}{\noindent}
\nc{\emp}{\emptyset}
\nc{\qed}{\hfill \lower.6em\hbox{\rule{1ex}{2ex}}}
\nc{\qued}[1]{\hfill \lower.6em\hbox{\rule{1ex}{2ex}{(#1)}}}
\nc{\ph}{\varphi}
\nc{\weg}{\scriptscriptstyle -}
\nc{\wegb}{\raisebox{.3ex}{$\,\scriptstyle -\,$}}
\nc{\gdw}{\ \Longleftrightarrow\ }
\nc{\gdwb}{\Longleftrightarrow}
\nc{\eps}{\hspace{-.2em}\in\hspace{-.2em}}
\nc{\epsb}{\in\hspace{-.2em}}
\nc{\dop}{:\!}
\nc{\cuep}{\stackrel{\cdot}{\cup}}
\nc{\neps}{\hspace{-.2em}\not\in\hspace{-.2em}}
\nc{\nepsb}{\not\in\hspace{-.2em}}
\nc{\spe}{\hspace{-.2em}\ni\hspace{-.2em}}
\nc{\ga}{\alpha} \nc{\gb}{\beta}\nc{\gd}{\delta}
\nc{\gl}{\lambda} \nc{\go}{\omega} \nc{\gs}{\sigma}
\nc{\gD}{\Delta}\nc{\gG}{\Gamma}
\font\msb=msbm10 scaled 1200
\font\msbs=msbm10 scaled 1000
\nc{\ess}{\mbox{\msb \char"28}}
\nc{\esss}{\mbox{\msbs \char"28}}
\nc{\sus}{\!\subset\!}
\nc{\sss}{\subset}
\nc{\oss}{\!\supset\!}
\nc{\gleich}{\! =\!}
\font\fract=eusb10 scaled 1200
\font\fracts=eusb10 scaled 1000
\nc{\kk}{\mbox{\fract K}} \nc{\FF}{\mbox{\fracts F}}
\nc{\LL}{\mbox{\fract L}} \nc{\DD}{{\sf D}}
\nc{\AS}{\mbox{$A$}} \nc{\BB}{\mbox{$B$}} \nc{\CC}{\mbox{$C$}}
\nc{\VS}{\mbox{\fract V}} 
\nc{\US}{\mbox{\fract U}} \nc{\PS}{\mbox{\fract P}}
\nc{\USs}{\mbox{\fracts U}}
\nc{\tief}{\rule[-.5ex]{0em}{1ex}}   
\nc{\unt}[1]{\underline{\tief {#1}}} 
\nc{\dt}{\gD\mbox{-tp}}
\nc{\leer}{\rule{0em}{1ex}}
\nc{\klg}{\:\raisebox{.4ex}{$\scriptstyle \leq$}\;}
\nc{\kklg}{\:\raisebox{.4ex}{$\scriptscriptstyle \leq$}\;}
\nc{\grg}{\:\raisebox{.4ex}{$\scriptstyle \geq$}\;}
\nc{\equi}{{\scriptstyle\equiv}\,}
\nc{\equ}{{\scriptstyle\equiv}\,}
\nc{\betw}{1\klg i\klg m} 
\nc{\JSL}{The Journal of Symbolic Logic}
\nc{\Def}{\abs \noindent{\bf Definition:} }
\nc{\Defi}{\abs {\bf Definition:} }
\nc{\Nota}{\abs \noi {\bf Notation:} }
\nc{\Notat}{\abs {\bf Notation:}J}
\nc{\Pro}{\noi {\bf Proof:}\ }
\nc{\Proo}{{\bf Proof:}\ }
\nc{\Rema}{\abs\noi{\bf Remark:}}
\nc{\Rem}{\abs {\bf Remark:}}
\nc{\Dann}{\Rightarrow}
\nc{\dann}{\rightarrow}
\nc{\wdann}{\rightarrow_{\mbox{$w$}}}
\nc{\card}{\raisebox{.4ex}{\scriptsize{\#}}}
\nc{\krmp}{1\klg r\klg m+1}
\nc{\kjr}{1\klg j\klg r}
\nc{\kin}{1\klg i\klg n}
\nc{\kkm}{1\klg k\klg m}
\nc{\kjm}{1\klg j\klg m}
\nc{\kklm}{1\klg k<l\klg m}
\centerline{\Large Many $\go$-categorical Structures Have}
\centerline{\Large the Small Index Property}
\abs
\centerline{BERNHARD HERWIG}
\centerline{Paris, France}
\centerline{June 1994}

\begin{quote}
{\bf Abstract}
\nl
Theorem: Let \AS\ be a finite $\kk_m$-free graph, 
$p_1,\ldots,p_n$  partial isomorphisms of \AS. Then
there exists a finite extension \BB, which is also  a $\kk_m$-free
graph, 
and automorphisms $f_i$ of \BB\
extending the $p_i$'s. This theorem can be used to prove the small
index property for the generic countable graph of this
class. The same method also works for a certain 
class of 
continuum many non isomorphic $\go$-categorical countable 
digraphs.
\end{quote}
Hrushovski proved this theorem for the class of all finite
graphs [Hr]; the proof presented here is an extension of his
proof.
\Nota
Let $p$ be a partial mapping on a set $A$.
By $D(p)$ we denote the domain of $p$
by $R(p)$ the range of $p$ (so $D(p)\sus A$ and
$R(p)\sus A$). In this paper the partial mappings under
consideration will always be injective.
The edge relation will always be
called $R$; in the first part of the paper we will deal with
graphs, so $R$ will be symmetric and irreflexive, in the second
part we will handle digraphs, so $R$ will be antisymmetric and
irreflexive.

\Def
If \AS\ is a graph and $a$ is a point (possibly element of a
graph extending $A$) we denote by $N_A(a)$ the set of
neighbours of $a$, sometimes considered as a pure set and
sometimes considered as a subgraph of $A$: $N_A(a):=\{b\eps
A\mid aRb\}$. A graph $A$ is called $\kk_m$-free (for
$m\eps\go$), if $\kk_m$, which is
the complete graph with $m$ vertices, is not
embeddable into $A$, i.e.\ there does not exist
$a_1,\ldots,a_m\eps A$ such that $a_iRa_j$ (for $\kklm$).

\begin{Satz}
Let $m,n\grg 1$, let 
$A$ be a finite $\kk_m$-free graph. Let $p_1,\ldots,p_n$ be
partial isomorphisms on $A$. There exists a finite $\kk_m$-free
graph $B,\ B\oss A$, and $f_1,\ldots,f_n\eps\mbox{\rm Aut}(B)$,
such that $f_i\oss p_i$.
\end{Satz}

Before going through the formal proof, which looks a little bit
technical, let us describe the main ideas of the proof: 
\nl
We know that the theorem holds for $\kk_3$-free graphs [Hg]. The
proof in the general case will be by induction on $m$. Let us
introduce for every $a\eps A$ a new colour (i.e. a unary
predicate) $U_a$ such that for $b\eps A$ $b$ is of colour $U_a$
iff $bRa$. Now $A$ is $\kk_m$-free if and only if  you can not
embed $\kk_{m-1}$ ``uni-coloured'' into $A$, that means there does
not exist a colour $U_a$ and elements $a_1,\ldots, a_{m-1}$ in
$A$ such that $a_kRa_l\ (\kklm)$ and all the $a_i$'s are of
colour $U_a$. Thus one can reduce $\kk_m$-freeness conditions
to certain $\kk_{m-1}$-freeness conditions, if one works with
coloured graphs. Here are the main problems, one has to
overcome doing this reduction.
\begin{enumerate}
\item
With respect to the colours $p_i$ is not longer a partial
isomorphism, it is a ``permorphism'', i.e.\ it respects the
colours only up to a permutation $\chi_i$ of the colours
$((U_a)^{\chi_i}=U_{a^{p_i}})$. But as can be seen in [Hg]
Hrushovski's original proof of the theorem in the case of the
graphs [Hr] also works for permorphisms.
\item
$\chi_i$ is not  yet really a permutation of the set of colours
$\{U_a\mid a\eps A\}$, it is only a partial function. As in
[Hg] one overcomes this problem by doing a type realizing step
to get a nice graph $C\oss A$, afterwards one looks at the colours
$\{U_d\mid d\eps C\}$ and extends the partially already
defined functions $\chi_i$ to permutations of this set.
\item
If one extends the graph $A$ considered as a (uni-coloured
$\kk_{m-1})$-free graph to a (uni-coloured $\kk_{m-1}$)-free
graph $B$ and the $p_i$'s to $f_i$'s, how can one ensure that
$B$ is $\kk_m$-free? Take into account that $B$ is $\{U_d\mid
d\eps C\}$-coloured, so we don't have for every $b\eps B$ a
colour $U_b$ such that all neighbours of $b$  have
colour $U_b$. This problem disappears by miracle: The
resulting graph $B$ looks locally like $A$. Especially for
every $a\eps A$ the neighbours of $a$ in $B$ will still all
have colour $U_a$. Now any orbit of the automorphism group of
$B$ (as graph) will have an element inside $A$. So to check
$\kk_m$-freeness of $B$, one has only to look for copies of
$\kk_m$ having one element $a$ inside $A$, but such a copy
would lead to a copy of $\kk_{m-1}$ of colour $U_a$.
\item
To get a proper induction, one has to prove the theorem for
coloured graphs and permorphisms. Starting with a coloured
graph, we have to introduce a new set of colours $\{U_d\mid
d\eps C\}$. But (uni-coloured $\kk_m$)-freeness does not
exactly mean (uni-coloured $\kk_{m-1}$)-freeness with respect
to the new colours. It means precisely: there does not exist an old
colour $U$ and a new colour $U_a$ such that $a$ is of colour
$U$ and a $U_a$-coloured copy 
of $\kk_{m-1}$ which is at the same
time $U$-coloured. We will call such a combination $(U,U_a)$ a
critical colouring, and we have to avoid copies of $\kk_{m-1}$
which are critical coloured. This last problem and the
notational complication arising from the fact that we are
dealing with permorphisms rather then isomorphisms
will make the prove look rather technical.
\end{enumerate}
The graph $C\oss A$ we will get by the type realizing step will
have as additional feature that every element in $A$ has
exactly the same number of neighbours in $C$, i.e. the same
number of colours. In the proof we will maintain this
condition, even if this is not really necessary. If one erases
in the proof all statements saying something like ``$\card\US^j
(a)=d_j$'' one get a slightly shorter proof. 
\abs

Now the definitions which follow, and the version of the
theorem (i.e. Lemma~2) which will be provable by induction
should be sufficiently motivated.

\Def
Let $S$ be a relational language. Let $\chi$ be a permutation
of the symbols in $S$ mapping every symbol to a symbol with the
same arity. Let $A$ be a $S$-structure and $p$ be a partial
mapping on $A$. $p$ is called a $\chi$-permorphism, if for
every $r\eps\go$ and every $r$-ary relation $R$ in $S$ and
every $a_1,\ldots,a_r\eps D(p)$: $Ra_1\ldots a_r\iff R^\chi
a_1^p\ldots a_r^p$.

\Def 
Let $\US^1,\ldots,\US^r$ be a family of disjoint  finite sets 
 of unary predicates 
(called colours) and $d_1,\ldots,d_r$ be constants. Let
$\US:=\bigcup_{1\kklg j\kklg r}\US^j$, $\LL:=\US\cup\{R\}$. If
$A$ is a $\LL$-structure and $a\eps A$, $V\eps \US$, then we
write 
$a\eps V$ to indicate that the unary predicate $V$ (or rather its
interpretation in $A$) is true for $a$; we also write ``$a$ is of
colour $V$''.
\begin{enumerate}
\item
A \US-graph is a $\{R\}\cup\US$-structure \AS\ such that 
\AS\ considered as a $\{R\}$-structure is a graph and
furthermore for every
$a\eps\AS$ and every $j$ ($\kjr$):  
$\card\US^j(a)=d_j$. Here
$\US^j(a):=\{V\eps\US^j\mid a\eps V\}$ is the set of colours of
$a$. $\US (a):=\{ V\eps\US\mid a\eps V\}$.

\item 
Let $\US_c\sus\US^1\times\ldots\times\US^r$. $\US_c$ will be
called the set of critical colourings. We call $A$
$\US_c$-$\kk_m$-free (for $m\eps\go$), if
there does not exist a colouring
$(V_1,\ldots,V_r)\eps\US_c$ and elements $a_1,\ldots,
a_m\eps A$, such that $a_k R a_l$ (for $\kklm$) and
$a_k\eps V_j$ (for $\kkm,\ \kjr$).
\end{enumerate}
\abs

\begin{Lemma}
Let $m\grg 1$. Let $\US^1,\ldots,\US^r$ be disjoint sets of
colours (where $r\grg 0$). Let
$\US:=\bigcup_{1\kklg j\kklg r}\US^j$. Let
$\chi_i^j\eps\mbox{Sym}(\US^j)$ (for $\kin,\ \kjr$),
$\LL:=\US\cup\{R\}$
$\chi_i:=\bigcup_{0\kklg j\kklg r}\chi_i^j\eps\mbox{Sym}
(\LL)$, where $\chi_i^0$ is the
identity on $\{R\}$. Let furthermore
$\US_c\sus\US^1\times\ldots\times\US^r$ be  a set called
critical colourings. Suppose $\US_c$ is $\chi_i$-invariant (for
$\kin$). \abs
Let A be a $\US$-graph. We
suppose that for every $\kjr$ there is given a constant $d_j$,
such that for every $a\eps A$ $\card (\US^j(a))\
=d_j$. Suppose $A$ is
$\US_c$-$\kk_m$-free.
\abs
Let $p_1,\ldots, p_n$ be partial mappings on $A$ such that $p_i$ is
a $\chi_i$-permorphism.
We suppose further that for $a\eps A$ and $1\klg j\klg r$ 
there exists a colour
$U_a^j\eps\US^j$ such that for $b\eps A$: $b\eps U_a^j\gdw aRb$
and such that for $a\eps D_i=D(p_i)$:
$(U_a^j)^{\chi_i}=U_{a^{p_i}}^j$.
Let us finally suppose that for every $a,b\eps A$, if $a\not=
b$ then there exist an $i$ ($\kin$) such that
$a^{p_i}=b$. (This is not really necessary.)

\abs\abs 
Then there exist a $\US$-graph $B$,
$B\oss A$, $B$ $\US_c$-$\kk_m$-free, $f_1,\ldots,f_n\eps
\mbox{Aut}(B),$\nl 
$f_i\oss p_i$, $f_i$ a $\chi_i$-permorphism (for
$1\klg i\klg n$). $B$ satisfies in addition: $\forall b\eps
B\forall a\eps A:\ aRb\Dann\forall j(\kjr)\ b\eps U_a^j$.
Furthermore $B$ will be chosen to satisfy in addition:
$\forall b\eps B\exists f\eps\;<f_1,\ldots,f_n>:\
b^f\eps A$.
\end{Lemma}

\abs\noi
From  the Lemma follows the theorem:\nl
Let in the Lemma $r=0$, that means we just talk about
(uncoloured) graphs. $\US_1\times\ldots\times\US_r$ just
contains the empty tuple $\gl$ and we let $\US_c=\{\gl\}$, then
$\US_c$-$\kk_n$-freeness just means $\kk_n$-freeness. $\chi_i$
is the identity on $\{R\}$ and $\chi_i$-permorphism just means
isomorphism of graphs. We can suppose w.l.o.g.\ that for every
$a,b\eps A$ there exists an $i\ (\kin)$ such that $a^{p_i}=b$
for example by assuming that $\{p_1,\ldots,p_n\}$ equals the
set of all partial isomorphisms on $A$.
Thus the lemma for $r=0$ just
implies the theorem.

\abs\noi
Proof of the Lemma:\nl
The proof goes by induction on $m$. Lets first treat the case
$m=1$. \nl
The proof follows in this case the lines of the original
proof of Hrushovski [Hr]. The first step will be a type
realizing step.
A subset of $A$ determines a ``type over $A$'' (considered as a
pure graph) in
our context. 
\abs

\noi Claim:\nl
a) There exist a finite graph $C$, $C\oss A$ ($C$ is a pure
graph so $C\oss A$ means $A$ is a substructure of $C$ as 
graphs)  and a constant 
$c_0$ such that  for every $A_0\sus A$: $\card (\{c\eps C\mid
N_A (c)=A_0\})=c_0$.
\nl
b) There exist bijections $h_1,\ldots,h_n\eps\mbox{Sym} (C)$,
$h_i\oss p_i$, such that for every $a\eps D_i,\ b\eps C:\
aRb\gdw a^{p_i}Rb^{h_i}$.
\abs

\noi Proof of a):\nl
Let for $A_0\sus A$ $c_{A_0}:=\card\{c\eps A\mid N_A (c)=A_0\}$
and $c_0:= \mbox{max}\{c_{A_0}\mid A_0\sus A\}$. Now for
every $A_0\sus A$ add $( c_0-c_{A_0})$ many points which have as
set of neighbours $A_0$ to get $C$.
\abs

\noi
Proof of b):\nl
Fix $i$ $(\kin)$
We first check, that for every $D_0\sus D_i\  
(D_i:=D(p_i),\ R_i:=R(p_i))$\nl
$\card (\{c\eps C\mid N_{D_i}
(c)=D_0\})= \card (\{c\eps
C\mid N_{R_i} (c)=D_0^{p_i}\}$:
\[
\card (\{c\eps C\mid N_{D_i}
(c)=D_0\})
=\sum_{B\sus (A-D_i)}\card (\{c\eps C
\mid N_A (c)=D_0\cup B\})=
\]\vspace{-.5ex}\[
c_0\cdot 2^{\card (A-D_i)}= c_0\cdot 2^{\card
(A-R_i)}= 
\card (\{c\eps C\mid N_{R_i}
(c)=D_0^{p_i}\}).
\]

Furthermore, because $p_i$ is a partial
isomorphism of graphs,  $p_i$ maps 
$\{c\eps C\mid N_{D_i}
(c)=D_0\}\cap D_i$ to 
$\{c\eps C\mid N_{R_i}
(c)=D_0^{p_i}\}\cap R_i$. So we find $h_i\oss p_i$ $h_i$ a bijection of $C$
mapping for every $D_0\sus D_i$ 
$\{c\eps C\mid N_{D_i}
(c)=D_0\}$ to  $\{c\eps C\mid N_{R_i}
(c)=D_0^{p_i}\}$, but this means exactly that $h_i$ has the
property we want.\qued{Claim}
\abs\abs

Now we do a duplicator step:\nl
We fix $C$ and $h_1,\ldots, h_n$, which we get from the claim.
We let $\Gamma\sus\mbox{Sym}(\LL)\times \mbox{Sym} (C)$ be the
subgroup generated by the elements $\gamma_i= (\chi_i,h_i)\
(\kin)$. Note that if $\gamma= (\chi,h)\eps\Gamma$ then 
$\chi$ fixes $\US^j$ (for $\kjr$) and $\chi$ fixes $\US_c$. On
$A\times \Gamma $, we define the equivalence relation $\equiv$ to
be the symmetric, reflexive and transitive closure of
$E=\{((a^{p_i},\gamma),(a,\gamma_i\gamma)) \mid\kin,\ a\eps
D_i, \gamma\eps\Gamma\}$. \nl
We note some basic facts:
\begin{enumerate}
\item
If $(a,(\chi,h))\equ (a_s,(\chi_s,h_s))$ then $a^h=a_s^{h_s}$
and $(\US (a))^\chi= (\US (a_s))^{\chi_s}$
\item
If $(a,(\chi,h))\equ (a_s,(\chi_s,h_s))$ and
$(b,(\chi,h))\equ (b_s,(\chi_s,h_s))$ then \nl
$aRb\gdw a_sRb_s$
\item
If $(a,(\chi,h))\equ (a_s,(\chi_s,h_s))$ and $c\eps C$ then\nl
$aRc^{h^{-1}}\gdw a_sRc^{h_s^{-1}}$ and $( U^j_a)^{\chi}=
(U^j_{a_s})^{\chi_s}$. 
\end{enumerate}
Proof of the facts:
\begin{enumerate}
\item
It suffices to prove 1.\ in the case
$((a,(\chi,h)),(a_s,(\chi_s,h_s)))$ is actually in $E$, so
$a_s\eps D_i$ and $a_s^{p_i}=a$ 
and $(\chi_s,h_s)= (\chi_i,h_i)\cdot
(\chi,h)$. But then $a^h=a_s^{p_ih}=a_s^{h_ih}=a_s^{h_s}$ and
$\US (a)^\chi=\US (a_s^{p_i})^\chi=\US (a_s)^{\chi_i\chi}=\US
(a_s)^{\chi_s}$ because $p_i$ is a $\chi_i$-permorphism.
\setcounter{enumi}{2}
\item Again we can suppose that $a_s\eps D_i,\
a=a_s^{p_i}$ and $(\chi_s,h_s)= (\chi_i,h_i)
(\chi,h)$. So $a_sRc^{h_s^{-1}}\gdw a_sRc^{h^{-1}h_i^{-1}}\gdw
a_s^{p_i}Rc^{h^{-1}}$ by the condition (in claim b)) on $h_i$ 
and $(U_{a_s}^j)^{\chi_s}=
(U_{a_s}^j)^{\chi_i\chi}= (U_{a_s^{p_i}}^j)^\chi= (U_a^j)^\chi$
by the conditions on $U_a^j$ in the hypothesis of the lemma.
\setcounter{enumi}{1}
\item
$aRb\gdw aRb^{hh^{-1}}\gdw a_sRb^{hh_s^{-1}}\gdw
a_sRb_s^{h_sh_s^{-1}}\gdw a_sRb_s$,\nl by 1. and 3.
\end{enumerate}

Now we are ready to 
define a $\LL$-structure on $A\times\Gamma/\equ$:\nl
For $e,f\eps A\times\Gamma/\equ:\nl
eRf\gdw\exists\gamma\eps\Gamma\exists a,b\eps A\ e=
(a,\gamma)/\equ, \ f= (b,\gamma)/\equ$ and $aRb$.
\nl
For
$V\eps\US$ and $e\eps A\times\Gamma/\equ$ we define 
\nl
$e\eps V\gdw\exists (\chi,h)\eps\Gamma\;\exists a\eps A\ e=
(a,(\chi,h))/\equ$ and $a\eps V^{\chi^{-1}}$. 
\nl
We note that
\begin{enumerate}
\setcounter{enumi}{3}
\item
For $\gamma\eps\Gamma$ $a,b\eps A$: $( (a,\gamma)/\equ) R
( (b,\gamma)/\equ)\gdw aRb$.
\item
For $(\chi,h)\eps\Gamma$, $a\eps A,\ V\eps\US$:
$( (a,(\chi,h))/\equ)\eps V^\chi\gdw a\eps V$
\end{enumerate}
Proofs:
\begin{enumerate}
\setcounter{enumi}{3}
\item
follows directly from 2.
\item
follows from 1.:\nl
$(a,(\chi,h))/\equ\eps V^\chi\gdw \nl (\exists a_s,\chi_s,h_s\
(a,(\chi,h))\equ (a_s,(\chi_s,h_s))$ and $a_s\eps
V^{\chi\chi_s^{-1}})\gdw a\eps V$,\nl 
namely: $a_s\eps
V^{\chi\chi_s^{-1}}\gdw V^{\chi\chi_s^{-1}}\eps\US (a_s)\gdw 
V^\chi\eps\US
(a_s)^{\chi_s}\gdw 
\nl
V^\chi\eps\US (a)^\chi
\gdw V\eps\US
(a)\gdw a\eps V$. \qued{facts}
\end{enumerate}
We define a map $i:A\dann A\times\Gamma/\equ$ by $i (a):=
(a,1)/\equ$ (where $1$ is the unit element in the group
$\Gamma$). 
$i$ is injective: if $(a,1)\equ
(b,1)$ then by 1.\ $a=b$. By 4.\ and 5.\ $i$ is an
embedding of $A$ into $A\times\Gamma/\equ$ as
$\LL$-structures. By identifying we suppose $A\sus
A\times\Gamma/\equ$ and we set $B=A\times\Gamma/\equ$. $B$ is
$\US_c$-$\kk_1$-free. For suppose there  exists
$(a,(\chi,h))/\equ\eps B$ and $(V_1,\ldots,V_r)\eps\US_c$ such
that $(a,(\chi,h))/\equ\eps V_j$ (for $\kjr$), then by 5.\
$a\eps V_j^{\chi^{-1}}$ and
$(V_1^{\chi^{-1}},\ldots,V_r^{\chi^{-1}})\eps\US_c$, because
$\chi^{-1}$ fixes $\US_c$. But this contradicts
$\US_c$-$\kk_1$-freeness of $A$.\abs
We define the automorphism $f_i$ (for $\kin$) by
$((a,\gamma)/\equ)^{f_i}:= (a,\gamma\gamma_i)/\equ$. This is an
$\chi_i$-permorphism, because for $V\eps\US$ and
$(a,(\chi,h))/\equ\eps B$: 
\[
(a,(\chi,h))/\equ\eps V\gdw a\eps V^{\chi^{-1}}\gdw a\eps
(V^{\chi_i})^{\chi_i^{-1}\chi^{-1}}\gdw
\]\[
(a,(\chi\chi_i,hh_i))/\equ\eps V^{\chi_i}\gdw
((a,(\chi,h))/\equ)^{f_i}\eps V^{\chi_i}
\]
The mapping $\ph:\Gamma\dann\;<f_1,\ldots,f_n>\ $(sending
$\gamma_i$ to $f_i$), which is given by
$((a,\gamma')/\equ)^{\ph (\gamma)}= (a,\gamma'\gamma)/\equ$ is
a surjective homomorphism of groups, and for
$((a,\gamma)/\equ)\eps B$ we have $((a,\gamma)/\equ)^{\ph
(\gamma^{-1})}= (a,1)/\equ\eps A$.

Finally let $a\eps A,\ b\eps B, \ aRb$. So let $a=
(a,1)/\equ= (a_s,(\chi_s,h_s))/\equ$ and $b=
(b_s,(\chi_s,h_s))/\equ$. Because $a_sRb_s$ we have
$b_s\eps\US_{a_s}^j$ (in $A$). So by 5.\ $b\eps
(\US_{a_s}^j)^{\chi_s}$ but by 3.\
$(\US_{a_s}^j)^{\chi_s}=\US_a^j$. This finishes the case $m=1$.
\qued{$m=1$}
\abs\abs
Now we do the step of induction $m\dann m+1\ (m\grg 1)$:\nl
We have the set of colours $\US^1,\ldots,\US^r$ and we know
that $A$ is $\US_c$-$\kk_{m+1}$-free. By a type realizing step
and by introducing new colours we want to consider $A$ as
satisfying a certain $\kk_m$-freeness condition and then we
want to apply the Lemma for $m$.
\abs
A subset $A_0\sus A$ and a colouring $\US_0\sus\US$ determines
a type over $A$ in this context. But not all of the types are
realizable in $\US_c$-$\kk_{m+1}$-free graphs. We call a tuple
$(A_0,\US_0)$ realisable, if for all $j$ $(\kjr)$ $\card
(\US_0\cap\US^j)=d_j$ and if there does not exist
$(V_1,\ldots,V_r)\eps\US_c\cap\US_0\times\ldots\times\US_0$
and elements $a_1,\ldots,a_m\eps A_0$ such that $a_kRa_l$ (for
$\kklm$) and $a_k\eps V_j$ (for $\kkm,\
\kjr$).  Here are facts about realizability:

\begin{enumerate}
\item
If $B_0\sus A_0\sus A$ and if $(A_0,\US_0)$ is realisable, then
$(B_0,\US_0)$ is realisable.
\item
If $C\oss A$ ($C$ a $\US$-graph) is
$\US_c$-$\kk_{m+1}$-free, then for every $c\eps C$ \nl $(N_A
(c),\US (c))$ is realisable. In particular
for every $a\eps A$ $( N_A (a),\US (a))$ is
realisable. 
\item
If
$(\kin)$ and $D_0\sus D_i$ and $\US_0\sus\US$:
$(D_0,\US_0)$ is realisable iff $(D_0^{p_i},\US_0^{\chi_i})$ is
realisable.
\end{enumerate}
Proof of 3.:\nl
Suppose $(D_0,\US_0)$ is realisable. Now we have for all $j,\
(\kjr):\ \card (\US_0^{\chi_i}\cap\US^j)=\card
(\US_0\cap\US^j)^{\chi_i}=d_j$. Suppose there exists
$(V_1^{\chi_i},\ldots,V_r^{\chi_i})\eps\US_c\cap\US_0^{\chi_i}\times
\ldots\times\US_0^{\chi_i}$ and
$a_1^{p_i},\ldots,a_m^{p_i}\eps D_0^{p_i}$ such that
$a_k^{p_i}Ra_l^{p_i}$ (for $\kklm$) and 
$a_k^{p_i}\eps V_j^{\chi_i}$ (for $\kkm$ and
$\kjr$). But then
$(V_1,\ldots,V_r)\eps\US_c\cap\US_0\times\ldots\times\US_0$ and
$a_1,\ldots,a_m\eps D_0$ and $a_kRa_l$ and $a_k\eps V_j$, because
$\US_c$ is $\chi_i$-invariant and because $p_i$ is a
$\chi_i$-permorphism. This is in contradiction to the
realizability of $(D_0,\US_0)$, therefore
$(D_0^{p_i},\US_0^{\chi_i})$ must be realisable.
\qued{facts}
\abs\abs
Now we do the type realizing step:\nl
Claim:\nl
a) There exists a $\US_c$-$\kk_{m+1}$-free $\US$-graph
$C\oss A$ and for every $t\ (0\klg t\klg\card A)$ a constant
$c_t$ such that for every $A_0\sus A$ and every $\US_0\sus\US$:
\[
\card\{c\eps C\mid N_A (a)\oss A_0,\ \US (a)=\US_0\}=\left\{
\begin{array}{ll}
c_{\card A_0} & \mbox{if $( A_0,\US_0)$ is realisable} \\
0   & \mbox{otherwise}
\end{array}\right.
\]
b) There exist bijections $h_1,\ldots,h_n\eps\mbox{Sym} (C),\
h_i\oss p_i$ such that for every $V\eps\US$, for every
$b\eps C$, 
for every $i\ (\kin)$: $b\eps V\gdw b^{h_i}\eps V^{\chi_i}$
and such that 
for every $a\eps D_i$, $b\eps C\ aRb\gdw a^{p_i}Rb^{h_i}$.
\abs\noi
Proof of a):\nl
Let $T=\card A$. We construct graphs $A=C_T\sus
C_{T-1}\sus\ldots\sus C_0=C$ and constants $c_T,\ldots c_0$
such that for every $t \ (T\grg t\grg 0) $ and for every
$A_0\sus A$ with $\card A_0\grg t$ and for $\US_0\sus\US$: $
\card\{c\eps C_t\mid N_A (a)\oss A_0,\ \US (a)=\US_0\}=
c_{\card A_0}$ if $( A_0,\US_0)$ is realisable and such that
$C_t$ is $\US_c$-$\kk_{m+1}$-free. We set $c_T=0$. If
$c_r,C_r$ are already constructed (for $T\grg r\grg t$) and
$t\grg 1$ then we will construct $C_{t-1}$ by adding points
which have exactly $t-1$ neighbours, all of them in $A$: For
$A_0\sus A$
with $\card A_0=t-1$, and for $\US_0\sus\US$ such
that $(A_0,\US_0)$ is realisable we define
$c_{A_0,\USs_0}:=\card\{c\eps C_t\mid N_A (c)\oss A_0,\ \US
(c)=\US_0\}$ and we define
$c_{t-1}$ to be the maximum of all these $c_{A_0,\USs_0}$. Now
to get $C_{t-1}$
we add for every realisable $(A_0,\US_0)$ (with $\card A_0=t-1$)
$c_{t-1}-c_{A_0,\USs_0}$ many points, which have as set of
neighbours exactly $A_0$ and as set of colours $\US_0$. 
$C_{t-1} $ is a $\US$-graph and $\US_c$-$\kk_{m+1}$-free and
for every $A_0\sus A$ (with $t-1\klg \card A_0$), for
every $\US_0\sus\US$: $\card \{c\eps C_{t-1}\mid N_A (c)\oss A_0,\
\US (c)=\US_0\}=c_{\card A_0}$ if $(A_0,\US_0)$ is
realisable. This is true, because for $A_0\grg t$ we did not
change the set in question (by going from $C_t$ to $C_{t-1}$)
and if $\card A_0=t-1$ then $\card \{c\eps C_{t-1}\mid N_A
(c)\oss A_0,\
\US (c)=\US_0\}=
\card \{c\eps C_t\mid N_A (c)\oss A_0,\
\US (c)=\US_0\}+ (c_{t-1}-c_{A_0,\USs_0})=c_{A_0,\USs_0}+
(c_{t-1}-c_{A_0,\USs_0})=c_{t-1}$.
\abs
The proof of b) is similar as the proof of claim b) in the case
$m=1$.
\nl
Here it is crucial to check, that for every $D_0\sus D_i=D
(p_i)$ and every $\US_0\sus\US$ ($R_i:=R (p_i)$):
\[
\card\{c\eps C\mid N_{D_i} (c)\gleich D_0,\ \US (c)
\gleich\US_0\}=
\card\{c\eps C\mid N_{R_i} (c)\gleich D_0^{p_i},\ 
\US (c)\gleich\US_0^{\chi_i}\}
\]
This is done by downwards induction on the size of $D_0$. The
step of induction is in the case $(D_0,\US_0)$ is realisable
(otherwise both sets are empty; here we are using fact 3.:
$(D_0,\US_0)$ is realisable$\gdw (D_0^{p_i},\US_0^{\chi_i})$ is
realisable) as
follows: 
\begin{eqnarray*}
 & &\card\{c\eps C\mid N_{D_i} (c)\gleich D_0,\ \US (c)\gleich
\US_0\}\\
 &=&
\card\{c\eps C\mid N_A (c)\oss D_0,\US (c)\gleich\US_0\}\\
 & & 
- \sum_{D_0\esss E\sus D_i}
\card\{c\eps C\mid N_{D_i} (c)\gleich E,\US (c)\gleich\US_0\}
\\
 &=&
c_{\card D_0}-
\sum_{D_0\esss E\sus D_i}
\card\{c\eps C\mid N_{R_i} (c)\gleich E^{p_i},\ \US (c)\gleich
\US_0^{\chi_i}\}
\\
 &=&
\card\{c\eps C\mid N_A (c)\oss D_0^{p_i},\ \US (c)
\gleich\US_0^{\chi_i}\}
\\
 & & - 
\sum_{D_0^{p_i}\esss E'\sus D_i}
\card\{c\eps C\mid N_{R_i} (c)\gleich E',\ \US (c)
\gleich\US_0^{\chi_i}\}
\\
 &=&
\card\{c\eps C\mid N_{R_i} (c)\gleich D_0^{p_i},\ 
\US (c)\gleich\US_0^{\chi_i}\}
\end{eqnarray*}
\hspace{-.5 ex}
\qued{claim}
\abs\abs
Now we introduce a new set of colours:
$\US^{r+1}=\{U_d^{r+1}\mid d\eps C\}$, where $U_d^{r+1}$ is a
new unary predicate for every $d\eps C$. We define
$\chi_i^{r+1}\eps\mbox{Sym} (\US^{r+1})$ (for $\kin$) by
$(U_d^{r+1})^{\chi_i^{r+1}}:=U_{d^{h_i}}^{r+1}$.
We let $\US':=\US\cup\US^{r+1}$ and
$\chi_i':=\chi_i\cup\chi_i^{r+1}\eps\mbox{Sym}
(\US'\cup\{R\})$. $\US_c'\sus\US^1\times\ldots\times\US^{r+1}$
is defined by 
\nl$(V_1,\ldots,V_r,U_d^{r+1})\eps\US_c'\gdw
(V_1,\ldots,V_r)\eps\US_c$ and $d\eps V_j$ ($\kjr$) 
(in $C$). 
\nl
$\US'_c$ is $\chi'_i$-invariant, because
$\US_c$ is $\chi_i$-invariant, because of the definition of
$\chi_i^{r+1}$ and because of the property of $h_i$ in claim b)
($b\eps V\gdw b^{h_i}\eps V^{\chi_i}$).
\abs

The colours in $\US^{r+1}$ are
in a natural way interpreted in $A$: for $a\eps A$
and $d\eps C$ we define $a\eps U_d^{r+1}\gdw dRa$ (in
$C$). Now $A$ is a $\US'$-graph. 
We have for every $a\eps A$: 
\nl $\card (\US^{r+1}
(a))=\card (\{d\eps C\mid a\eps\US_d^{r+1}\})=\card
(\{d\eps C\mid aRd\})$

$=\sum_{\USs_0\sus\USs}\card\{d\eps C\mid N_A
(d)\oss\{a\},\ \US (d)=\US_0\}=k_a\cdot c_1$.
\nl
Here $c_1$ is a constant appearing in claim a), and
\nl
$k_a=\card\{\US_0\sus\US\mid (\{a\},\US_0) \mbox{is
realisable}\}$. But for $a,b\eps A$ $k_a=k_b=:k$: pick $i$ such
that $a^{p_i}=b$, then 
$\{\US_0\sus\US\mid (\{a\},\US_0) \mbox{is
realisable}\}^{\chi_i}=
\{\US_0\sus\US\mid (\{b\},\US_0) \mbox{is
realisable}\}
$. (This argument is only needed in the case $m+1=2$.) 
Thus we define the constant $d_{r+1}:=k\cdot c_1$. 
\abs

$A$ is $\US'_c$-$\kk_m$-free. Otherwise there would exist
$a_1,\ldots a_m\eps A$ and
\nl $(V_1,\ldots,V_r,U^{r+1}_d)\eps\US'_c$ 
such that $a_kRa_l$ (for
$\kklm$) and $a_k\eps V_j$ (for $\kkm$, $\kjr$) and
$a_k\eps U^{r+1}_d$. But then $(V_1,\ldots,V_r)\eps\US_c$,
$d\eps V_j$, $a_kRd$. This means $a_1,\ldots,a_m,d$ contradicts
the $\US_c$-$\kk_{m+1}$-freeness of $C$.

\abs
For $1\klg i\klg n$ $p_i$ is a
$\chi_i'$-permorphism. Furthermore 
the conditions on the colours
$\US^{r+1}_a$ (for $a\eps A$) are satisfied, e.g. for
$a\eps D_i$ we have
$(U_a^{r+1})^{\chi_i}=U_{a^{h_i}}^{r+1}=U^{r+1}_{a^{p_i}}$ 
because $p_i\sus h_i$. 
\abs 

By the Lemma for $m$ we find a finite $\US'$-graph
$B$, $B\oss A$, $B$ $\US'_c$-$\kk_m$-free and $f_1,\ldots,
f_n\eps\mbox{Aut} (B)$, $f_i\oss p_i$, $f_i$ a
$\chi_i'$-permorphism 
having the indicated properties. Now we consider $B$ just as a
$\US$-graph. The only thing we still have to check, is that $B$
is $\US_c$-$\kk_{m+1}$-free.
\abs

Suppose there exist $(V_1,\ldots,V_r)\eps \US_c$ and elements
$a_0,\ldots,a_m\eps B$ such that $a_kRa_l$ and $a_k\eps
V_j$. W.l.o.g\ we suppose $a:=a_0\eps A$: Otherwise choose
$f\eps\;<f_1,\ldots,f_n>$, such that $a_0^f\eps A$ and choose
$\chi'\eps\mbox{Sym} (\LL')$ such that $f$ is a
$\chi'$-permorphism and let $\chi=\chi'\mid_{\LL}$; now still
$(V_1^\chi,\ldots,V_r^\chi)\eps\US_c$ and $a_k^fRa_l^f$ and
$a_k^f\eps V_j^\chi$. Because $a\eps V_j$ we  have
$(V_1,\ldots,V_r,U_a^{r+1})\eps\US_c'$. The additional
condition on $B$ implies now: $a_k\eps\US_a^{r+1}\
(\kkm)$. Thus we get a contradiction to the
$\US_c'$-$\kk_m$-freeness of $B$.
\qued{Lemma,Theorem} 

\begin{Theorem}
Let $m\grg 1$, let $M$ be the generic countable $\kk_m$-free
graph. $M$ has the Small Index Property
\end{Theorem}
Proof: The proof is given in [HHLS]. The use there of Hrushovski's
Lemma in the proof of the Small Index Property for the generic
graph must be replaced by theorem~1.
\qued{Theorem}

\abs\abs
Now we are turning to the case of the digraphs:\nl
Lets suppose there is a fixed (possibly infinite) family $\FF$ 
of finite tournaments
(i.e. of digraphs F such that for every $a,b\eps F$
if $a\not=b$ then $aRb$ or $bRa$). Lets look at the class $\kk$
of all finite $\FF$-free digraphs $A$, i.e. of digraphs 
such that no $F\eps\FF$ is embeddable into
$A$. This class has the (free) Amalgamation Property, and the
resulting generic countable digraph $M_{\FF}$ will be called
a Henson digraph. Henson [Hen] showed that there continuum many
non isomorphic such digraphs. 
\begin{Theorem}
a) Let $A\eps \kk$, let $p_1,\ldots,p_n$ be partial
isomorphisms on $A$, then there exist $B\eps\kk$, $B\oss A$ and
$f_1,\ldots,f_n$ automorphisms of $B$, $f_i\oss p_i$.
\nl
b) The Small Index Property holds for $M_{\FF}$. I.e.\/ it holds
for all Henson digraphs.
\end{Theorem}
Proof:
\nl
Again The proof given in [HHLS] shows that b) follows from
a). 
\nl
It suffices to prove a) in the case that \FF\ is a finite
class:  Let \FF\ be infinite and let $A\eps\kk$. Let $m:=\card
A$. Let $\FF_1$ be a finite family of tournaments of size $m+1$
containing every isomorphism
type of such a tournament. 
We define $\FF_0:=\{F\eps\FF\mid \card F\klg m\}\cup\FF_1$ and
we assume w.l.o.g $\FF_0$ to be finite. Now $A$ is $\FF_0$-free
and every $\FF_0$-free graph $B$ is also \FF-free, i.e. in
$\kk$.

Now we want to prove a) in the case $\FF$ is finite by
induction on the maximal size of a tournament in \FF. This
prove will be very similar to the prove of theorem~1. We only
want to point out the differences.
Again the following lemma
will be the
``permorphism version'' of the Theorem, which is provable by
induction.
Here we will have no
restriction on the cardinality of colours a single point can
have.
\begin{Lemma}
Let \US\ be a finite set of colours, let $\chi_1,\ldots,\chi_n\eps
\mbox{Sym} (\US)$. Let $T_1,\ldots,T_m$ be finite tournaments,
$T_j=\{t^j_1,\ldots,t_{l_j}^j\}$, let for $\kjm$
$\US_j\sus(\mbox{Pot}
(\US))^{l_j},\ \US_j\
\chi_i$-invariant. Let $A$ be a $\US$-coloured digraph such
that for every $j$ ($\kjm$) $A$ is  $\US_j$-$T_j$-free
(i.e.\ there does not exist an embedding of digraphs
$i:T_j\dann A$, $s_k:=i (t^j_k)$ such that $(\US
(s_1),\ldots,\US (s_{l_j}))\eps\US_j$). Let $p_1,\ldots,p_n$ be 
partial mappings such that $p_i$ is a $\chi_i$-permorphism. 
\abs
Then there exists a finite $\US$-coloured digraph $B$, $B\oss A$, 
$B$ $\US_j$-$T_j$-free ($\kjm$) and 
$f_1,\ldots,f_n$ total permorphisms on $B$,
$f_i\oss p_i$.
\abs\noi
$B$ satisfies in addition:
For every mapping $U:A\dann \US$ (we will write $U_a$ instead of
$U(a))$ \nl
if for $b\eps A: \ (bRa\gdwb b\eps U_a)$ (rsp.: $aRb
\gdwb b\eps U_a$) 

and if for $a\eps D_i=D(p_i)$: 
$(U_a)^{\chi_i}=U_{a^{p_i}}$ \nl
then for every $b\eps B$, for every $a\eps A$:
$(bRa\Dann b\eps U_a)$ (rsp.: $aRb\Dann b\eps
U_a$).
\end{Lemma}
Proof:
\nl
The proof of the lemma in this version of the preprint
will be rather sketchy. The proof goes by induction on the maximal
size of the tournaments $T_j$ ($\kjm$) involved. In the case
that the 
maximal size is 1 it goes exactly like the proof  
in the case of the 
graphs (case $m=1$). So let's suppose that the maximal size is
$>1$. 
\abs
Let us first introduce some notation. If $A$ is a digraph and
$a$ is a point (possibly element of a digraph extending $A$) we
denote by $N^+_A (a):=\{b\eps A\mid aRb\}$ and by $N^-_A
(a):=\{b\eps A\mid bRa\}$.

 Let $T_1,\ldots,T_r$ be the tournaments of size $>1$ and 
$T_{r+1},\ldots,T_m$ be the tournaments of size $1$. We write for
$\kjr$ $T_j^+:=N^+_{T_j} (t^j_1)$ and $T_j^-:=N^-_{T_j}
(t^j_1)$, thus $T_j-\{t^j_1\} \gleich T^+_j\cup T^-_j$.

We will call $C$ (where $C\oss A$)
free of critical copies of $T_j$ if there is no
embedding  $\rho:\ T_j 
\dann C$ such that (if one denotes $s_l=\rho(t_l^j)$) 
$s_2,\ldots,s_{l_j}\eps A$ and such that $(\US(s_1),\ldots,\US(s_{l_j}))
\eps\US_j$.

\abs
Now the first step will be, to find $C\oss A$, such that $C$ is
free of critical copies of $T_j$ and such  that there are bijections 
$h_1,\ldots,h_n\eps\mbox{Sym}(C)$  such that for every $V\eps\US$, 
for every  $c\eps C$: $c\eps V\gdw c^{h_i}\eps V^{\chi_i}$ and 
such that for every $a\eps D_i,\ c\eps C:\ (aRc\gdw a^{p_i}Rc^{h_i})$ 
and $(cRa\gdw c^{h_i}R a^{p_i})$.

In this case a subset $A_0^+\sus A$ and a disjoint subset
$A_0^-\sus A$ and a set of colours $\US_0$ determine a type
over $A$. We call $(A_0^+,A_0^-,\US_0)$ realisable, if there
does not exist $j\klg r$ and an embedding (as digraphs) 
$i:T_j^+\cup T_j^-\dann A_0^+\cup A_0^-$ such that $i
(T_j^+)\sus A_0^+$ and $i(T_j^-)\sus A_0^-$ and such that
$(\US_0,\US (s_2),\ldots,\US(s_{l_j}))\eps\US_j$ (again
$s_l=i(t_l^j))$. Note that for $C\oss A$: $C$ is free of
critical copies of $T_j$ (for every $j\klg r$) 
iff for every $c\eps C$:
$(N_A^+(c),N_A^-(c),\US (c))$ is realizable.
Now again choose constants $c_t$ ($1\klg t\klg\card A$) and
choose $C$  to ensure that for every $A_0^+,\ A_0^-\sus A,\
\US_0\sus \US$:
\[
\card (\{c\eps C\mid N_A^+(c)\oss
A_0^+,\ N_A^-(c)\oss
A_0^-,\ \US(c)=\US_0 \})
\]\vspace{-.5ex}\[=
\left\{
\begin{array}{ll}
c_{\card ( A_0^+\cup A_0^-)} & \mbox{if $( A_0^+,A_0^-,\US_0)$ 
is realisable} \\
0   & \mbox{otherwise}
\end{array}\right.
\]
As in the proof of Lemma~2 one proves that bijections
$h_1,\ldots,h_n$ with the desired properties exist.

\abs

Now our new set of colours will be $\US'=\US\cuep\{U^+_c\mid
c\eps C\} \cuep\{U^-_c\mid
c\eps C\}$ (with new colours $U^+_c,U^-_c$). For $a\eps A$ we
define ($a\eps U_c^+$ iff $cRa$) and ($a\eps U_c^-$ iff
$aRc$. This way we equip $A$ with a $\US'$-structure and $p_i$
is a $\chi_i'$-permorphism, if we define
$(U_c^+)^{\chi_i'}\gleich U^+_{c^{h_i}}$ and
$(U_c^-)^{\chi_i'}\gleich U^-_{c^{h_i}}$
 (for $\kin$). 
For $1\klg j\klg r$ we set $T'_j=\{t_2^j,\ldots,t_{l_j}^j\}$
and we say (for $\VS_2,\ldots,\VS_{l_j}\sus\US'$): 
$(\VS_2,\ldots,\VS_{l_j})\eps\US'_j$ iff there exists $c\eps
C$ such that for $2\klg l\klg l_j$ $U_c^{\epsilon_l}
\eps\VS_l$ and
such that $(\US(c),\VS_2\cap\US,\ldots,\VS_{l_j}
\cap\US)\eps\US_j$, where $U_c^{\epsilon_l}=\US_c^+$, if $t^j_l\eps
T^+_j$ and $\US_c^{\epsilon_l}=\US_c^-$, if $t^j_l\eps
T^-_j$. 

$A$ is $\US_j'$-$T_j'$-free.
\abs

By induction there is an $\US'_j$-$T_j'$-free (for $\kjr$)
digraph $B$, which is also $\US_j$-$T_j$-free (for $r<j\klg m$),
which satisfies all we want. The only thing which still needs
checking is that $B$ is $\US_j$-$T_j$-free (for $\kjr$). 

Lets suppose there exists a $j$ $\kjr$ and there exists
an embedding $i:T_j\dann B$, such that
$(\US(s_1),\ldots,\US(s_{l_j}))\eps \US_j$ (if we write
$s_l:=t^j_l$). We can suppose that $s_1\eps A$. We consider
(for $\epsilon=+$ and for $\epsilon=-$) the
mapping $U^\epsilon:A\dann \US'$ such that $U^\epsilon
(a)=U^\epsilon_a$, and we fix $c=s_1\eps A$.
Take $l,\ 2\klg l\klg l_j$; if $t_l^j\eps T_j^+$, then
$t_1^jRt_l^j$, so $cRs_l$, so by the additional property $B$
satisfies: $s_l\eps U_c^+$. The same way one checks that if
$t_l^j\eps T_j^-$, then $s_l\eps U_c^-$. Now we consider the
embedding $i$ restricted to $T_j'=T_j-\{t_1^j\}$.
Now it is easy to check
that $(\US' (s_2),\ldots,\US' (s_{l_j}))\eps\US_j'$. But this
is in contradiction to the 
$\US'_j$-$T_j'$-freeness of $B$.\qued{Lemma,Theorem}
\abs\abs\abs

Finally we want to mention that there is a more general
theorem stating the possibility of extending partial
isomorphisms. It is a theorem about relational structures in any
relational language. To state it we need some notation.
Let $S$ be a finite relational language.
\Seite

\Def
\begin{itemize}
\item
Let $L$ be a $S$-structure. $L$ is called a link structure if
it consists of just one element or if there is a $k$-ary
($k\eps\go$) relation symbol $R$ in $S$ and a $k$-tuple
$(a_1,\ldots,a_k)$ in $L$ such that $Ra_1\ldots a_k$ and such
that $a_1,\ldots,a_k$ contains all elements of $L$. 
\item
Let $A$ and $B$ be $S$-structures. $A$ and $B$ have the same
link type if for every link structure $L$: $L$ is embeddable
into $A$ iff $L$ is embeddable into $B$.
\item
Let $T$, $A$ be $S$-structures, let $\rho\dop T\dann A$ be a
function. $\rho$ is called a weak homomorphism (notation: $\rho\dop
T\wdann A$) if for every $R$ in $S$ ($R$ $k$-ary) and every
$s_1,\ldots, s_k\eps T$: If $Rs_1\ldots s_k$ (in $T$) then
$Rs_1^\rho\ldots s_k^\rho$ (in $A$).
\item
Let $\FF$ be a  set of finite $S$-structures. Let $A$ be
a $S$-structure. $A$ is called $\FF$-free if there does not
exist $T\eps\FF$ and $\rho\dop T\wdann A$.
\end {itemize}
For example if $S$ consists of just one binary relation symbol,
and $A$ is a graph (rsp. digraph) and if the $S$-structure 
$B$ has the same link
structure as $A$ then $B$ is a graph (rsp. digraph); but this is
not true for tournaments.
\begin{Satz}
Let $\FF$ be a finite set of finite $S$-structures. Let $A$ be
a finite $\FF$-free $S$-structure. Let $p_1,\ldots,p_n$ be
partial isomorphisms on $A$. There exist a finite $\FF$-free
$S$-structure $B$, $A\sus B$ and automorphisms
$f_1,\ldots,f_n$ on $B$ ($f_i\oss p_i$) such that $B$ and $A$ 
have the same link type.
\end{Satz}
The proof of the theorem is just the translation of the proofs
of theorem 1 and theorem 4 a) into this more general
context. Note that a weak homomorphisms mapping $\kk_m$ to a
graph $A$ (and a weak homomorphisms mapping a tournament into a
digraph) is necessarily  an embedding.

\Seite
\subsection*{References}
[C] P. Cameron, Oligomorphic Permutation Groups, LMSLNS 152,
\nl\indent Cambridge University Press, 1990,\bbs\noi
[Hg] B.Herwig, Extending Partial Isomorphisms, Combinatorica,
to appear,\bbs\noi
[Hr] E. Hrushovski, Extending Partial Isomorphisms of Graphs,
\nl\indent Combinatorica 12(1992), 411-416\bbs\noi  
[HHLS] Hodges, Hodkinson, Lascar, Shelah: The Small Index
Property for 
\nl\indent $\go$-stable, $\go$-categorical, structures and for the
random graph,\nl\indent Journal of the LMS 48(1993), 204-218. \bbs\noi
[T] J.K. Truss, Generic Automorphisms of Homogeneous
Structures,\nl\indent preprint, 1990.

\end{document}